\documentclass[12pt]{amsart}

\newtheorem{theorem}{Theorem}[section]
\newtheorem*{theorem A}{Theorem A}

\theoremstyle{remark}
\newtheorem{remark}{Remark}[section]
\theoremstyle{remark}

\theoremstyle{definition}

\newtheorem{definition}{Definition}[section]

\numberwithin{equation}{section}
\def\({\left ( }
\def\){\right )}
\def\<{\left < }
\def\>{\right >}
\def\e{\eqref}

\begin{document}

\title[Classification of spherical Lagrangian submanifolds]{Classification of spherical Lagrangian submanifolds in complex Euclidean spaces}

\author[B.-Y. Chen]{Bang-Yen Chen}

  \address{Department of Mathematics, Michigan State University, 619 Red Cedar Road, East Lansing, Michigan 48824--1027,  U.S.A.}
\email{bychen@math.msu.edu}

\begin{abstract}   An isometric  immersion  $f:M^n\to \tilde M^n$ from a Riemannian $n$-manifold $M^n$ into a K\"ahler $n$-manifold $\tilde M^n$ is called {\it Lagrangian\/} if the  complex structure $J$ of the ambient manifold $\tilde M^n$ interchanges each tangent space of $M^n$ with the corresponding normal space. In this paper, we completely classify spherical Lagrangian submanifolds in complex Euclidean spaces. In this paper, we also provide two corresponding classification theorems for Lagrangian submanifolds in the complex pseudo-Euclidean spaces with arbitrary complex index.
\end{abstract}


\maketitle

\section{Introduction.}

\def\ii{\hskip.013in {\rm i}}
\def\e{\eqref}
  
A symplectic manifold is a smooth manifold, $M$, equipped with a closed non-degenerate differential 2-form, $\omega$, called the symplectic form. Isotropic submanifolds  (or totally real submanifolds in  \cite{CO}) are submanifolds where the symplectic form restricts to zero, i.e. each tangent space is an isotropic subspace of the ambient manifold's tangent space. The most important case of the isotropic submanifolds is that of Lagrangian submanifolds. A Lagrangian submanifold is, by definition, an isotropic submanifold of maximal dimension, namely half the dimension of the ambient symplectic manifold.

Lagrangian submanifolds arise naturally in many physical and geometric situations.  One important example is that the graph of a symplectomorphism in the product symplectic manifold $(M\times M,\omega\times -\omega)$ is Lagrangian.
By Darboux's theorem and its generalizations, it is well known that the extensions of a $k$-manifold $M$ to a $2n$-symplectic manifold in which
$M$ is isotropic are classified, up to a local symplectomorphism about $M$, by the isomorphism
classes of $2(n-k)$-dimensional symplectic vector bundles over $M$ (see,
for instance, \cite[page 24]{W} or \cite{CM}). This result implies that there are no obstructions to local Lagrangian immersions without further restrictions.
 On the other hand, the author obtains obstructions of Lagrangian isometric immersions
from Riemannian $n$-manifolds into complex space forms in terms of his $\delta$-invariants (cf. \cite{1993,2000,book}).
Consequently, an important problem in the theory of Lagrangian
submanifolds in K\"ahler geometry is to determine Lagrangian submanifolds  with ``some given special  geometric properties''. For instance, a method
was given in \cite{J1}  to construct an important
family of special Lagrangian submanifolds in {\bf C}$^n$ with large symmetric groups.  In \cite{a1,a2}, R. Aiyama introduced a spinor-like representation formula which  parameterizes  immersions through two complex functions $F_1,F_2$ and a real one (the Lagrangian angle $\beta$). Aiyama's formula is  useful to construct many examples of  Lagrangian surfaces in {\bf C}$^2$. Also, the author has completely determined all Lagrangian surfaces of constant curvature in complex space forms in a series of articles (\cite{c2004}--\cite{c7}).

One fundamental question concerning Lagrangian submanifolds of a complex Euclidean $n$-space ${\bf C}^n$ is the following:
\vskip.05in 

{\it ``Which Lagrangian submanifolds of ${\bf C}^n$ are spherical\,?''}

\vskip.05in 

In this paper, we obtain a complete solution to this fundamental question. More precisely, we completely classify all spherical Lagrangian submanifolds  of  complex Euclidean spaces. 

Our main result is the following.

\begin{theorem} \label{T:1} Let $L:M^n\to {\bf C}^n$ be a Lagrangian isometric immersion of a Riemannian $n$-manifold $M^n$ into the complex Euclidean $n$-space ${\bf C}^n$.  Then $L(M^n)$ is contained in a hypersphere of ${\bf C}^n$ if and only if the following two statements hold:
\vskip.05in

{\rm (1)} $M^n$ is an open portion of the Riemannian product of a circle $S^1$ and a Riemannian $(n-1)$-manifold $N^{n-1}$. 
\vskip.05in 

{\rm (2)} Up to dilations and translations, the immersion is given by
\begin{align} L(t,u_2,\ldots,u_n)=e^{\ii t} \phi(u_2,\ldots,u_n),\end{align}
where $\phi$ is the horizontal lift of a Lagrangian isometric immersion $\hat \phi:N^{n-1}\to CP^{n-1}(4)$ of a Riemannian $(n-1)$-manifold $N^{n-1}$ into the complex projective space $CP^{n-1}(4)$ of constant holomorphic sectional curvature $4$.
\end{theorem}

In the last section, we provide two corresponding classification theorems for Lagrangian submanifolds in complex pseudo-Euclidean spaces with arbitrary complex index.

\section{Preliminaries.} We follow the notations from \cite{c1973,book}.

\subsection{Basic formulas and definitions}
Let $f:M\to {\bf C}^n$ be an isometric immersion of a Riemannian $n$-manifold $M$ into the complex Euclidean $n$-space  ${\bf C}^n$. We denote the Riemannian connections of $M$ and  ${\bf C}^n$ by $\nabla$ and $\tilde \nabla$, respectively; and by $D$  the connection on the normal bundle.

The formulas of Gauss and Weingarten are given respectively by
\begin{align}  \label{2.1} &\tilde\nabla_X Y =\nabla_X Y + h(X,Y),\\&\tilde\nabla_X \xi =- A_\xi X + D_X \xi\end{align}
  for tangent vector fields $X$ and $Y$ and normal  vector field  $\xi$. 
The second fundamental form $h$ is related to the shape operator $A_\xi$ by
\begin{align}\<h(X,Y),\xi\> = \<A_\xi X,Y\>.\end{align}

If we denote the Riemann curvature tensor of $\nabla$  by $R$, then the equations of Gauss and Codazzi are given respectively by
\begin{align} &\label{2.3}\<R(X,Y)Z,W\> =  \< h( X,W),h(Y,Z)\> -
\< h( X,Z),h(Y,W)\>,\\&  (\nabla h)(X,Y,Z) = (\nabla
h)(Y,X,Z),\end{align}
where $X,Y,Z,W$ are vector fields tangent to $M$ and $\nabla h$ is defined by 
\begin{align}  \notag (\nabla h)(X,Y,Z) = D_X h(Y,Z) - h(\nabla_X
Y,Z) - h(Y,\nabla_X Z).\end{align}

When $M$ is a Lagrangian submanifold in {\bf C}$^n$,  we also have
(cf. \cite{CO})
\begin{align} \label{2.7} &D_X JY = J \nabla_X Y,\\&
\<h(X,Y),JZ\>=\<h(Y,Z),JX\>=\<h(Z,X),JY\>.
\end{align}

\subsection{Horizontal lift of Lagrangian submanifolds}

The following link between Legendrian submanifolds and Lagrangian submanifolds is due to  \cite{R} (see also \cite[pp. 247--248]{book}).
\vskip.1in

 {\bf Case} (i): $CP^n(4)$.  Consider Hopf's fibration $\pi :S^{2n+1}\to CP^n(4).$  For a given point $u\in S^{2n+1}(1)$, the horizontal space at $u$ is the orthogonal complement of $\i u, \, \i=\sqrt{-1},$ with respect to the metric  on $S^{2n+1}$ induced from the metric on ${\bf C}^{n+1}$.   
Let $\iota : N \to CP^n(4)$ be a Lagrangian isometric immersion. Then there is a covering map
$\tau: \hat N \to N$ and a horizontal  immersion $\hat \iota :\hat N \to S^{2n+1}$ such that
$\iota\circ \tau=\pi \circ \hat \iota$.  Thus each Lagrangian immersion can be lifted locally (or globally if $N$ is simply-connected) to a Legendrian immersion of the same Riemannian manifold. In particular, a minimal Lagrangian submanifold of $CP^{n}(4)$ is lifted to a minimal Legendrian submanifold of the Sasakian $S^{2n+1}(1)$.

Conversely, suppose that $f: \hat N \to S^{2n+1}$ is a Legendrian isometric immersion. Then $\iota =\pi\circ f:  N\to  CP^n(4)$ is again a Lagrangian isometric immersion.  Under this correspondence the second fundamental forms $h^f$ and $h^\iota$ of $f$ and $\iota$ satisfy $\pi_*h^f=h^\iota$.  Moreover, $h^f$ is horizontal with respect to $\pi$. 

\vskip.1in

{\bf Case} (ii): $CH^n(-4)$. We consider  the complex number space  ${\bf C}^{n+1}_1$ equipped
with the pseudo-Euclidean metric:
$$g_0=-dz_1d\bar z_1 +\sum_{j=2}^{n+1}dz_jd\bar z_j .$$
 Consider 
$$H^{2n+1}_{1}(-1)=\{z\in {\bf C}^{2n+1}_1: \<z,z\>=-1\}$$ with the canonical Sasakian structure, where $\<\;\,,\;\>$ is the induced inner product. We consider  the complex number space  ${\bf C}^{n+1}_1$
Put \begin{align}&T'_z=\{ u\in {\bf C}^{n+1}:\<u,z\>=0\},\;\;\\& H_1^1=\{\lambda\in {\bf  C}: \lambda\bar\lambda=1\}.\end{align}
 Then there is an $H^1_1$-action on $H_1^{2n+1}(-1)$,  $z\mapsto \lambda z$ and at each point $z\in H^{2n+1}_1(-1)$, the vector $\xi=-\i z$ is tangent to the flow of the action. Since the metric $g_0$ is Hermitian, we have $\<\xi,\xi\>=-1$.  The quotient space $H^{2n+1}_1(-1)/\sim$, under the identification induced from the action, is the complex hyperbolic space $CH^n(-4)$ with constant holomorphic sectional curvature $-4$ whose complex structure $J$ is  induced from the complex structure $J$ on ${\bf C}^{n+1}_1$ via Hopf's fibration $\pi :H^{2n+1}_1(-1 )\to  CH^n(4c).$

 Just like case (i), suppose that $\iota: N \to CH^n(-4)$ is a Lagrangian immersion, then there is an isometric covering map
$\tau: \hat N \to N$ and a Legendrian immersion $f: \hat N \to H_1^{2n+1}(-1)$ such that $\iota\circ \tau=\pi\circ f$.  Thus every Lagrangian immersion into $CH^n(-4)$ an be lifted locally (or globally if $N$ is simply-connected) to a Legendrian immersion into $H^{2n+1}_1(-1)$. In particular,  Lagrangian minimal submanifolds of $CH^{n}(-4)$ are lifted to Legendrian minimal submanifolds of $H^{2n+1}_{1}(-1)$.
Conversely, if $f:\hat N \to H_1^{2n+1}(-1)$ is a Legendrian immersion, then $\iota =\pi\circ f:  N\to CH^n(-4)$ is a Lagrangian immersion.  Under this correspondence the second fundamental forms $h^f$ and $h^\iota$ are related by $\pi_*h^f=h^\iota$. Also, $h^f$ is horizontal with respect to $\pi$. 

Let $h$ be the second
fundamental form of $M$ in $S^{2n+1}(1)$ (or in $H^{2n+1}_1(-1)$). Since $S^{2n+1}(1)$ and $H^{2n+1}_1(-1)$ are totally umbilical  with one as its mean curvature in  ${\bf C}^{n+1}$ and in ${\bf C}^{n+1}_1$,
respectively, we have  
\begin{align}\label{2.10} \hat \nabla_XY=\nabla_XY+h(X,Y)-\varepsilon L,\end{align} 
where $\varepsilon=1$ if the ambient space is ${\bf C}^{n+1}$; and $\varepsilon=-1$ if it is ${\bf C}^{n+1}_1$.

\section{Proof of Theorem \ref{T:1}}
Let $L:M^n\to {\bf C}^n$ be a Lagrangian isometric immersion from a Riemannian $n$-manifold $M^n$ into ${\bf C}^n$. Denote by $J$ the complex structure of ${\bf C}^n$. Suppose that $L(M^n)$ lies in a hypersphere of ${\bf C}^n$. Then, by applying suitable dilation and translation, $L(M^n)$ lies in the unit hypersphere $S^{2n-1}(1)$ centered at the origin. Let ${\bf x}$ denote the position function of $M^n$ in ${\bf C}^n$ and let $V=J{\bf x}$. Then {\bf x} is a normal vector field of $M^n$ and $V$ is tangent to $M^n$.

If we put $$\mathcal D_1={\rm Span}\{V\},\;\; \mathcal D_2=(\mathcal D_1)^\perp,$$ then $TM^n={\mathcal D}_1\oplus {\mathcal D}_2$.
Since ${\bf C}^n$ is K\"ahlerian, we have 
\begin{align}\label{3.1}\tilde \nabla_Z V=J\tilde\nabla_Z {\bf x}=JZ\end{align}for any $Z\in TM^n$.
From \e{3.1} we find
\begin{align}\label{3.2} &\nabla_Z V=0,\\& \label{3.3}h(Z,V)=JZ,\;\;\\&\label{3.4} h(V,V)=JV=-{\bf x}, \end{align}
for any $ Z\in TM^n$.

It follows from \e{3.2} that $\mathcal D_1$ is a totally geodesic distribution, i.e., $\mathcal D_1$ is integrable and its leaves are totally geodesic in $M^n$. Moreover,  we also know from \e{3.4} that each integral curve of $V$ is an open portion of a great circle of $S^{2n-1}(1)$.

Since $\mathcal D_1$ and $\mathcal D_2$ are orthogonal complementary distributions, it follows from \e{3.2} that $\<[X,Y],V\>=0$.
Thus, $\mathcal D_2$ is an integral distribution.

From the orthogonality of $\mathcal D_2$ and $V=J{\bf x}$, we find
\begin{align}\label{3.5}\<\right.\!\tilde \nabla_XY, V\! \left.\>=-\<\right.\! Y, \tilde\nabla_XV\! \left.\>=-\<\right.\!Y, JX\! \left.\>=0.\end{align}
On the other hand, it follows from formula \e{2.1} of Gauss that
\begin{align}\label{3.6}\<\right.\!\tilde \nabla_XY, V\! \left.\>=\<\right.\! \nabla_XY, V\! \left.\>.\end{align}
Therefore, after combining \e{3.5} and \e{3.6}, we conclude that $\mathcal D_2$ is  a totally geodesic distribution. Consequently, according to the well-known de Rham decomposition theorem, $M^n$ is locally the Riemannian product $S^1\times N^{n-1}$ of a circle $S^1$ and a Riemannian $(n-1)$-manifold $N^{n-1}$. Hence, there exists a local coordinate system $\{t,u_2,\ldots,u_n\}$ on $M^n$ such that $V=\frac{\partial}{\partial t}$ and $\frac{\partial}{\partial u_2},\ldots,\frac{\partial}{\partial u_n}\in \mathcal D_2$.

It follows from \e{3.2}, \e{3.3} and \e{3.4} that the immersion $L:M^n\to {\bf C}^n$ satisfies the following system of partial differential equations:
\begin{align}\label{3.7} &\frac{\partial^2 L}{\partial t^2}=\ii \frac{\partial L}{\partial t},\\& \label{3.8}\frac{\partial^2 L}{\partial t\partial u_j}=\ii \frac{\partial L}{\partial u_j},\;\;\\&\label{3.9} \frac{\partial^2 L}{\partial u_i\partial u_j}=\sum_{k=2}^n \Gamma^k_{ij} \frac{\partial L}{\partial u_k}+h\(\frac{\partial }{\partial u_j},\frac{\partial }{\partial u_j}\)-\<\frac{\partial }{\partial u_j},\frac{\partial }{\partial u_j}\>L, \end{align}
for $i, j=2,\ldots, n$, where $\Gamma^k_{ij}$ are the Christoffel symbols.

After solving \e{3.7} we get
\begin{align}\label{3.10} &L=A(u_2,\ldots,u_n)+e^{\ii t}\phi(u_2,\ldots,u_n)\end{align}
for some ${\bf C}^n$-valued functions $A$ and $\phi$. 

By substituting \e{3.10} into \e{3.8} we find $$\frac{\partial A}{\partial u_j}=0,\;\; j=2,\ldots,n.$$ Thus, $A$ is a constant vector, say $c_0$. Hence, \e{3.10} becomes
\begin{align}\label{3.11} &L=c_0+e^{\ii t}\phi(u_2,\ldots,u_n).\end{align}

From \e{3.11} we get 
\begin{align}\label{3.12} &\frac{\partial L}{\partial t}=\ii e^{\ii t}\phi,\end{align}
which implies that $\<\phi,\phi\>=1$, since $V$ is a unit vector field. Therefore, it follows from \e{3.10} and $\<L,L\>=1$ that
\begin{align}\label{3.13} &0=\<c_0,c_0\>+2\<c_0,e^{\ii t}\phi\>.\end{align}
Now, after taking the differentiation of  \e{3.13} twice with resect to $t$, we obtain $\<c_0,e^{\ii t}\phi\>=0$. Therefore, we derive from \e{3.13} that $c_0=0$. Consequently, \e{3.11} implies that $L$ takes the form:  
\begin{align}\label{3.14} &L=e^{\ii t}\phi(u_2,\ldots,u_n).\end{align}

From \e{3.14} we have
\begin{align}\label{3.15} &\frac{\partial L}{\partial u_j}=e^{\ii t}\frac{\partial \phi}{\partial u_j},\;\; j=2,\ldots,n,\end{align}
which gives 
\begin{align}\label{3.16} &\<\frac{\partial L}{\partial u_i},\frac{\partial L}{\partial u_j}\>=\<\frac{\partial \phi}{\partial u_i},\frac{\partial \phi}{\partial u_j}\>,\;\; i,j=2,\ldots,n,\end{align}
Therefore, $\phi$ is an isometric immersion.

Next, after applying the Lagrangian condition of the immersion $L:M^n\to {\bf C}^n$, we may conclude from \e{3.15} that $\phi:N^{n-1}\to S^{2n-1}(1)$ is a Legendrian isometric immersion. Consequently, $\phi$ is a horizontal lift of a Lagrangian isometric immersion $\hat \phi:N^{n-1}\to CP^{n-1}(4)$ from $N^{n-1}$ into $CP^{n-1}(4)$.

The converse can be verified by direct computation. \qed

\section{Lagrangian submanifolds in ${\bf C}^n_s$}

Let $n$ and $s$ be integers satisfying $n\geq 1$ and $0\leq s\leq n$. The complex manifold ${\bf C}^n$ endowed with the real part of the
Hermitian form
\begin{align} b_{s,n}(z,w)=-\sum_{j=1}^s \bar z_jw_j+\sum_{j=s+1}^n\bar z_jw_j\;\;{\rm for}\; z,w\in {\bf C}^n\end{align}
defines a flat indefinite complex space form of index $2s$, denote by  ${\bf C}^n_s$.

For $c>0$, the pseudo hypersphere $S^{2n-1}_{2s}(c)$ in ${\bf C}^n_s$ is defined by
\begin{align} S^{2n-1}_{2s}(c)=\Big\{{\bf z}=(z_1,\ldots,z_n)\in {\bf C}^n_s: \<{\bf z},{\bf z}\>=\frac{1}{c}>0\Big\}.\end{align}
Similarly, for $c<0$, the pseudo hyperbolic space $H^{2n-1}_{2s-1}(c)$ in ${\bf C}^n_s$ is defined by
\begin{align} H^{2n-1}_{2s-1}(c)=\Big\{{\bf z}=(z_1,\ldots,z_n)\in {\bf C}^n_s: \<{\bf z},{\bf z}\>=\frac{1}{c}<0\Big\}.\end{align}

On $S^{2n-1}_{2s}(c)$ (resp., $H^{2n-1}_{2s-1}(c)$), we define an endomorphism $P$ of the tangent bundle of $S^{2n-1}_{2s}(c)$ (resp., of $H^{2n-1}_{2s-1}(c)$) by
\begin{align}  JX=PX+FX, \;\; X\in T(S^{2n-1}_{2s}(c))\;\;  \hbox{(resp.,  $X\in T(H^{2n-1}_{2s-1}(c))$},\end{align} where $J$ is the complex structure of ${\bf C}^n_s$ and $PX$ and $FX$ denote the tangential and normal components of $JX$, respectively. 

Let $V$ denote the tangent vector field of $S^{2n-1}_{2s}(c)$ (resp., of $H^{2n-1}_{2s-1}(c)$) given by $V=J{\bf x}$, where ${\bf x}$ is the position function of $S^{2n-1}_{2s}(c)$ (resp., the position function of $H^{2n-1}_{2s-1}(c)$) in ${\bf C}^n_s$.

Analogous to Legendrian submanifolds in $S^{2n-1}$, we make the following definition.

\begin{definition} An $(n-1)$-dimensional submanifold $N^{n-1}$ of $S^{2n-1}_{2s}(c)$ (resp., of $H^{2n-1}_{2s-1}(c)$) is called a {\it Legendrian submanifold}  if the following two conditions are satisfied:
\vskip.05in

(1) $V$ is normal to $N^{n-1}$ and

\vskip.05in
(2) $P$ maps each tangent vector of $N^{n-1}$ into a normal vector.
\end{definition}

\begin{remark} {\rm Every Legendrian submanifold $N^{n-1}$ of $S^{2n-1}_{2s}(c)$ (resp., of $H^{2n-1}_{2s-1}(c)$) has index $s$ (resp. has index $s-1$).}\end{remark}

For Lagrangian submanifolds in ${\bf C}^n_s$, we have the following two results,

\begin{theorem} \label{T:4.2}  Let $L:M_s^n\to {\bf C}^n_s$ be a Lagrangian isometric immersion of a pseudo-Riemannian $n$-manifold $M^n_s$ into the complex pseudo-Euclidean $n$-space ${\bf C}^n_s$. 
Then $L(M^n_s)$ is contained in the unit pseudo hypersphere $S^{2n-1}_{2s}(1)\subset {\bf C}^n_s$ if and only if the following two statements hold:

\vskip.05in

{\rm (a)} $M^n_s$ is the direct product of  a circle $S^1$ and an $(n-1)$-dimensional pseudo Riemannian manifold $N^{n-1}_s$.
\vskip.05in

{\rm (b)} The immersion is given by
\begin{align} L(t,u_2,\ldots,u_n)=e^{\ii t} \psi(u_2,\ldots,u_n),\end{align}
where $\psi: N^{n-1}_s\to S^{2n-1}_{2s}(1)\subset {\bf C}^n_s$ is a Legendrian isometric immersion from $N^{n-1}_s$ into $S^{2n-1}_{2s}(1)$.
\end{theorem}

\begin{theorem} \label{T:4.3}  Let $L:M_s^n\to {\bf C}^n_s$ be a Lagrangian isometric immersion of a pseudo-Riemannian $n$-manifold $M^n_s$ into the complex pseudo-Euclidean $n$-space ${\bf C}^n_s$. 
Then $L(M^n_s)$ is contained in the pseudo hyperbolic space $H^{2n-1}_{2s-1}(-1)\subset {\bf C}^n_s$ if and only if the following two statements hold:
\vskip.05in

{\rm (a)} $M^n_s$ is the direct product of a time-like line $\mathbb E^1_1$ and an $(n-1)$-dimensional pseudo Riemannian manifold $N^{n-1}_{s-1}$.
\vskip.05in

{\rm (b)} The immersion is given by
\begin{align} L(t,u_2,\ldots,u_n)=e^{\ii t} \rho(u_2,\ldots,u_n),\end{align}
where $\rho: N^{n-1}_{s-1}\to H^{2n-1}_{2s-1}(-1)\subset {\bf C}^n_s$ is a Legendrian isometric immersion from $N^{n-1}_{s-1}$ into $H^{2n-1}_{2s-1}(-1)$.
\end{theorem}

Since Theorems \ref{T:4.2} and \ref{T:4.3} can be proved in the same way as Theorem \ref{T:1} with only minor modification, so we omit their proofs.

\end{document}